\documentstyle[12pt]{article}
\input psfig
\font\bit=cmssi12 at 12truept
\newtheorem{thm}{Theorem}
\newtheorem{dfn}{Definition}
\newtheorem{lem}{Lemma}
\newtheorem{prop}{Proposition}

\newtheorem{clm}{Claim}
\newtheorem{fct}{Fact}
\newtheorem{rmk}{Remark}

\newcommand{\Z}{{\bf Z}}
\newcommand{\N}{{\bf N}}
\newcommand{\Hu}{\tilde{\cal H}}
\newcommand{\Hl}{{\cal H}}
\newcommand{\Rh}{{\bf R}^{(h)}}
\newcommand{\Rv}{{\bf R}^{(v)}}
\newcommand{\Sh}{{\bf S}^{(h)}}
\newcommand{\Sv}{{\bf S}^{(v)}}
\newcommand{\Ss}{{\bf S}^{(\sigma)}}
\newcommand{\T}{{\bf T}^2}
\newcommand {\superarrow}[2]{\begin{array} {c} #1 \\ \longrightarrow \\ #2 \end{array}}
\newcommand {\ceil}[1]{\lceil #1 \rceil}
\newcommand {\nin}{\not\in}

\newcommand{\R}{{\bf R}}
\begin{document}

\begin{center} {\large\bf A toral diffeomorphism with a nonpolygonal rotation
set.}\\ \bigskip \bigskip

 Jaroslaw Kwapisz \footnote{Supported in part by
Polish Academy of Sciences grant 210469101 ``Iteracje i Fraktale''.}\\
Mathematics Department\\ SUNY at Stony Brook\\ NY 11794-3651\\ \bigskip
December, 1992 \footnote{Revised in April 1995.} \end{center}

\bigskip \bigskip

\begin{abstract}
We construct a diffeomorphism of the  two-dimensional torus which is isotopic to the identity and whose rotation set is not a polygon. 
\end{abstract}

\bigskip

\begin{center} {\large\bf Introduction.}  \end{center}

\bigskip

The rotation set is an important invariant associated with a 
homeomorphism of the two-dimensional torus. In the literature
there are a number of different definitions of the rotation set, 
however all of these sets can only differ on the boundary.
Moreover, all of the sets have the same convex
hull which is equal to the Misiurewicz-Ziemian  rotation set. This rotation set consists of all asymptotic average winding vectors exhibited by the system and is always a compact convex subset of the plane (\cite{MZ-1}).

A natural question is whether every compact convex subset of the plane is a
rotation set for some torus homeomorphism. This question is still open.
Even deciding which line segments are rotation sets is not yet resolved.
Given any interval with rational slope passing through a point with rational coordinates it is easy to give a homeomorphism that has the interval as its rotation set. The construction in \cite{FM} implies that intervals with at least one rational endpoint and irrational slope can be realized as well.
Based on their results for time-one maps of toral flows Franks and Misiurewicz conjecture that these two are the only possibilities (\cite{FM}).

In the case of nonempty interior, it is known that all polygons with vertices having rational coordinates can be realized as rotation sets
(\cite{JK-1}).
When the rotation set has interior, it depends continuously on the map
(\cite{MZ-2}).
This clearly implies that some sets which are not rational polygons must occur,
however it is still possible that only polygons occur.
In this note we show that this is not the case.
We construct a $C^1$-diffeomorphism whose rotation set is not a polygon.
The regularity of the example is determined by the use of 
a circle homeomorphism exhibiting a wandering arc. By Denjoy's theorem
wandering arcs do not happen in $C^2$-category (\cite{Denjoy},
\cite{Nitecki}). One is prompted then to ask whether there is a 
$C^{\infty}$-diffeomorphism whose rotation set has infinitely
many extremal points?

\begin{figure}[htp]
\centerline{\psfig{figure=rotset.fig,width=0.7\hsize}}
\centerline{Figure 0.}
\end{figure}

The construction of a $C^1$-diffeomorphism, $G$, whose rotation set
is shown in Figure 0, goes roughly as follows.
Pick a irrational number $\rho$. 
The one-dimensional skeleton of the two dimensional
torus is a bouquet of two circles. Call one of them
the {\bit vertical circle} and the other the {\bit horizontal circle}.  
The definition of $G$ begins by defining a map $F$ on this bouquet.
The map $F$ is the composition of two transformations.
The first transformation is obtained by mapping the
horizontal circle into itself by 
a circle homeomorphism which has a rotation number $\rho$ and exhibits a wondering arc containing the common point of the two circles.
Since we want the transformation to be continuous we must
fold a piece of the vertical circle into the horizontal one, and
we keep this folding as simple as possible. 
The second transformation is defined analogously with the 
roles of the horizontal and vertical circles interchanged. 

The dynamics of the composition of the two
transformations, $F$, will contain trajectories which
remain on a fixed circle. These
orbits all have a common average rotation rate $\rho$ and give rotation vectors $(0,\rho)$ and $(\rho,0)$.
The trajectories which visit the vertical and
horizontal circles repeatedly account for all other extremal points
except $(0,0)$. Our ability to effectively calculate the rotation vectors of these trajectories 
stems from the fact that the essential part of the dynamics is
controlled by an infinite Markov partition. Indeed, for any given
natural number $n$, the set of points on, say, the  horizontal circle 
that jump to the vertical one after exactly $n$ iterates is an arc. 
The collection of thus defined arcs has the Markov property.
Section 2 contains the details of the construction of the map
$F$ on the bouquet of circles.

The next step is to ``perturb'' $F$ to a
$C^1$-smooth embedding of a neighborhood of the 
skeleton into itself. The perturbation is carefully chosen so
that its dynamics are as close as possible to those of the 
inverse limit of the map on the skeleton. 
An obvious difficulty arises at the trajectory of the
common point of the vertical and horizontal circles. This
is overcome by locking this point into a wandering domain.
This is where the construction depends crucially on
the use of a Denjoy example.  It is a routine
to conclude that the rotation set of the perturbed map is the
same as that of $F$. 
To complete the construction, we extend the dynamics to
$G$ on the whole torus by putting a source in the complement of the
neighborhood of the skeleton. This places the vector $(0,0)$
in the rotation set of $G$.
The construction of  the ``perturbed'' map with appropriate proofs can be
found in Section 3. For a more general approach yet yielding only continuous embedding see Section 4.

\bigskip
\begin{center}
{\large\bf Section 1: Preliminaries and statements of results.}
\end{center}
\bigskip

Let ${\bf T}^2$ be the 2-dimensional torus  ${\bf R}^2 / {\bf Z}^2$ and $\pi : {\bf R}^2 \rightarrow {\bf T}^2$ the projection. For any compact subset $X \subset {\bf T}^2$, let $\tilde{X} := \pi ^ {-1}(X)$ and denote by $\Hu (X)$  the set of all continuous  $ \tilde{F} : \tilde{X} \rightarrow  \tilde{X}$ satisfying  \( \tilde{F}(\tilde{x}+v)=\tilde{F}(\tilde{x})+v \) for all $\tilde{x} \in \tilde{X}, v \in {\bf Z}^2$.\
Denote by  $\Hl (X)$ all mappings  $F:  X \rightarrow X$ that are projections of those in  $\Hu (X)$. Notice that $\Hl ({\bf T}^2)$ is the space of all continuous maps of ${\bf T}^2$ homotopic to the identity. For any $\tilde{F} \in \Hu (X)$ we have the {\bit displacement function} $\phi_{\tilde{F}}: X \rightarrow {\bf R}^2$ defined by $\phi_{\tilde{F}}(x) := \tilde{F}(\tilde{x})-\tilde{x}$ where $\tilde{x}$ is a lift of $x$. 
We borrow the following definition from \cite{MZ-1}.

\begin{dfn} \hspace{-0.09in}{\bf .}
For $ \tilde{F} \in \Hu (X) $ and a point $x \in X$,  its {\bit rotation set} $\rho(\tilde{F}, x)$ is the set of all limit points of a sequence $(\frac{1}{n}(\tilde{F}^n (\tilde{x})-\tilde{x}))_{n \in \N}$ where $\tilde{x}$ is a lift of $x$. This set is independent on the choice of lift $\tilde{x}$.
The {\bit rotation set of the map} is  
\[ \rho ( \tilde{F}) := \bigcap_{m>0} cl\left( \bigcup_{n>m} \{\frac{1}{n}(\tilde{F}^n (\tilde{x})-\tilde{x}): \tilde{x} \in \tilde{X} \} \right) \subset {\bf R}^2. \]
For $F \in \Hl (X)$ we set $\rho(F) := \rho(\tilde{F})$ where $\tilde{F}$ is a lift of $F$. This defines $\rho(F)$ only up to translation by a vector in ${\bf Z}^2$.
\end{dfn}

It is a useful general fact that the convex hull of the rotation set for any map is generated by the rotation sets of generic points of ergodic invariant measures.
In particular, we have the following proposition.
(Here we write $Conv(A)$ for the convex hull of $A \subset {\bf R}^2$.) 

\begin{prop}[\cite{MZ-1}] \hspace{-0.09in}{\bf .}
For any $\tilde{F} \in \Hu (X)$ the following sets coincide :
\begin{itemize}
\item[(i)] $Conv(\rho(\tilde{F}))$;
\item[(ii)] $\{\int \phi_{\tilde{F}}  \ d\mu : \mu \ is \ an \ F-invariant\ probability \ measure \ on \ X\}$;
\item[(iii)] $Conv \left( \bigcup {\rho(\tilde{F},x) : x\in X \ such \ that \ \rho(\tilde{F},x) \ is \ a  \ point} \right) $.
\end{itemize}
\end{prop}
\bigskip

Although we defined the rotation set for maps on subsets of the torus our primary 
interest is in the case when $X={\bf T}^2$ and $\tilde{F} \in \Hu ({\bf T}^2)$ is a homeomorphism. In this setting $\rho (\tilde{F})$ is a convex compact in  ${\bf R}^2$(\cite{MZ-1}). Our main result is existence of a $C^1$-diffeomorphism  $ \tilde{G} \in \Hu ( {\bf T^2}) $ with $\rho (\tilde{G})$ which is not a polygon. To be more precise let us adopt the following definitions. (We use $\lceil x \rceil$ for the smallest integer greater or equal then $x$ and $\lfloor x \rfloor$ for the largest integer less or equal then  $x$.) 

\begin{dfn} \hspace{-0.09in}{\bf .}
For any irrational $\rho \in (0,1)$ we define a sequence of numbers $(\alpha_n)_{n \in \N}$ by the formula  
\[ 
 \alpha_n := \lceil n \rho \rceil - n \rho
\]
and vectors ${\bf \rho}_{m,n}$ in ${\bf R}^2$ by 
\[ 
{\bf \rho}_{m,n} := \frac{(\lceil m \rho \rceil, \lceil n \rho \rceil)}{m+n+1}.
\]
Also we define the convex set $\Lambda_{\rho}$ by
\[ \Lambda_{\rho} := Conv \{0, \; {\bf \rho}_{m,n}: m,n \in {\bf N} ,\   \alpha_m < \rho , \ \alpha_n  < \rho \}. \]
\end{dfn}

Figure 0 gives an indication of the shape of $\Lambda_{\rho}$.
The next proposition asserts a crucial property of this set, namely, that it is not a polygon.

\begin{prop} \hspace{-0.09in}{\bf .}
For any irrational $\rho \in (0,1)$, the set $\Lambda_{\rho}$ has infinitely many extremal points. There are exactly two accumulation points of the  extremal points of $\Lambda_{\rho}$, namely,  $(0, \rho)$ and $(\rho,0)$.
\end{prop}

\bigskip

{\bf Proof of Proposition 2.}
Recall that $ {\bf \rho}_{m,n} := (x_{m,n},y_{m,n})$ where 
\[
x_{m,n} :=  \lceil m \rho \rceil / (m+n+1) = (m \rho + \alpha_m)/(m+n+1),
\]
\[ 
y_{m,n}:=\lceil n \rho \rceil / (m+n+1) ) = (n \rho + \alpha_n)/(m+n+1).
\] 
Suppose  $m_k , n_k \in {\bf N}$, $k = 1,2,3 ... $ are such that $x_{m_k,n_k} \rightarrow 0$. Then $n_k/m_k \rightarrow \infty$, so $ y_{m_k, n_k}$ is asymptotic to  $ \lceil n_k \rho \rceil / n_k$ and thus converges to $\rho$.
On the other hand, sequences  $(m_k), (n_k) $  with ${\bf \rho}_{m_k,n_k} \in \Lambda_\rho$ for which $x_{m_k,n_k} \rightarrow 0$ exist due to irrationality of $\rho$. In this way the vectors $\rho_{m,n}$ in $\Lambda_{\rho}$ which approach the $y$-axis accumulate to $(0, \rho)$. Also, since $\Lambda_\rho$ is contained in the first quadrant $\{ (x,y) : x,y \geq 0 \} $ and contains the origin, it follows that $(0, \rho)$ is an extremal point. 

Denote by $\gamma_{m,n}$ the slope of the line through ${\bf \rho}_{m,n}$ and $( \rho, 0)$. We claim the following:

If $m_k, n_k $ with ${\bf \rho}_{m_k,n_k} \in \Lambda_\rho$ for $k =1,2,3...$, are such that 
\[ 
 \lim_{k \rightarrow \infty} \gamma_{m_k,n_k} = \sup \{  \gamma_{m,n} : {\bf \rho}_{m,n} \in \Lambda_\rho \}, 
 \]
then $\gamma_{m_k,n_k}$ is not eventually constant and  $  \lim_{k \rightarrow \infty} {\bf \rho}_{m_k,n_k} = (0, \rho)$.

To prove this claim, note that we have 
\[
 \gamma_{m_k,n_k} = \frac {n_k \rho + \alpha_{n_k} - \rho (m_k+n_k+1)}{
m_k \rho + \alpha_{m_k}} \ ,
\] 
\[
 \gamma_{m_k,n_k} = -1 + \frac{\alpha_{m_k} + \alpha_{n_k }- \rho}{m_k \rho + \alpha_{m_k}}.
\]
Since $\rho$ is irrational $ \alpha_{n_k}$'s can be arbitrarily close to $\rho$
so it is clear that the supremum is $> -1$ and is not attained. Thus either $ m_k \rightarrow \infty $ or $ n_k \rightarrow \infty $. The first possibility leads to $\gamma_{m_k,n_k} \rightarrow  -1$ which is definitely not the supremum. It follows that $m_k$'s are bounded and $n_k \rightarrow \infty$ which gives $ {\bf \rho}_{m_k,n_k} \rightarrow (0, \rho) $. 

This last fact guarantees that $ (0, \rho) $ is a condensation point of extremal points. The same is true of $(\rho,0)$. In this way we have infinite number of extremal points. 
Moreover, using the formula for the slope $\gamma_{m,n}$ one can see that points ${\bf \rho}_{m,n}$ accumulate towards the anti-diagonal $\{ (x,y) : x+y=\rho \}$. Clearly there are no other  extremal points on it  except for $(0,\rho)$ and   $(\rho,0)$.
This proves the second part of Proposition 2. $\Box$

Now we are ready to  formulate our theorem.

\begin{thm}[Main Theorem]  \hspace{-0.09in}{\bf .}
For any irrational $\rho \in (0,1)$ there is a $C^1-$diffeomorphism  $ G \in \Hl ({\bf T}^2)$ whose rotation set $\rho(G)$ equals $\Lambda_{\rho}$ $(mod \ {\bf Z}^2)$.  
\end{thm}

An outline of the proof of Theorem 1 can be found in the introduction.
Sections 2 and 3 give more formal arguments.

\bigskip
\begin{center}
{\large\bf Section 2: A map on a one-dimensional skeleton of the torus.}
\end{center}
\bigskip

Denote the $x$ and  $y$-axis in ${\bf R}^2$  by $\Rh$  and $\Rv$ respectively. These lines project to circles on  $\T$ which we call $\Sh$ and $\Sv$. We put  $X := \Sh \cup \Sv$ and $\tilde{X} := \pi^{-1} (X)$.
We have obvious injections $ {\bf R}/{\bf Z} = {\bf S}^1  \rightarrow \Ss$ for $ \sigma \in \{h,v\}$. For $ x \in {\bf S}^1 $ we denote the corresponding point in $\Ss$ by  $ x^{(\sigma)} \in  \Ss$, $\sigma \in \{h,v\} $. We also have a ``projection'' $X \rightarrow  {\bf S}^1$ sending $x \mapsto \underline{x}$ that is defined by the relation  $ \underline{x^{(\sigma)}} = x, \ x \in X$.

Define the set $\Omega_{\rho} \subset {\bf R}^2$ by :
\[
\Omega_\rho := Conv \{ {\bf \rho}_{m,n} : \alpha_m, \alpha_n < \rho \},
\]
where ${\bf \rho}_{m,n}, \ \alpha_m, \ \alpha_n $ are as defined in Section 1.

We construct $F: X \rightarrow X$  , $F \in \Hl(X)$ with a rotation set 
\[ \rho (F) = \Omega_\rho \ (mod \ {\bf Z}^2).\]

This map is a composition of two analogous transformations of X. Roughly speaking, the first one is what one gets trying to rotate $\Sh$ by a positive angle keeping the map continuous and the unavoidable folding of $\Sv$ into $\Sh$ as simple as possible. The other does the same to $\Sv$. The detailed definition preceded by some necessary preliminaries follows.

 We think of the one-dimensional circle ${\bf S}^1$ as ${\bf R}  \ (mod \ {\bf Z})$.
Let $I$ be a small symmetric closed arc around zero and $\tilde{I}$ its lift to $\R$ containing zero.  For $\phi : {\bf S}^1 \rightarrow {\bf S}^1 $ a degree-one map that is 1-1 everywhere except $I$ where it has a plateau, let $ \tilde{\phi} : {\bf R} \rightarrow {\bf R} $ be the lift satisfying $ \tau := \tilde{\phi}(\tilde{I}) \in [0,1]$, (see Figure 1). As long as we are not concerned with the issue of smoothness, to make our construction work we need to pick any such a $\phi$ with the rotation number $\rho( \tilde{\phi})$ equal to $\rho$.

\begin{figure}[htp]
\centerline{\psfig{figure=degree-one.fig,width=0.7\hsize}}
\centerline{Figure 1.}
\end{figure}

However, since we ultimately want a $C^1$-example  we have to choose $\phi$ more carefully. We may start with a Denjoy $C^1$-diffeomorphism of ${\bf S}^1$ with rotation number equal to $\rho$. We may further require that it has a wandering domain  slightly larger  than $I$, say equal to the dilation of I by a factor of $18/17$ about the origin  that we denote by $18/17 \cdot I$. Now since $18/17 \cdot I$ is wandering, any modification of  the diffeomorphism on this arc without altering its image does not affect the rotation number. In particular, we can redefine the map on $18/17 \cdot I$ so that  $I$ is sent to a point and $C^1$-smoothness is preserved. This is the map we are ultimately  going to take for $\phi$. Let us stress at this point that we will not use wandering properties of $\phi$ until the considerations of Section 2.

There is a unique homeomorphism $\psi : {\bf S}^1 \rightarrow {\bf S}^1$ such that $\phi = \psi \circ p$ where $p$ collapses $I$ to zero and is affine on the complement of $I$. Indeed, $\psi(x) := \phi \circ p^{-1}(x)$ for $x \neq 0$ extends through $0$ because $\phi \circ p^{-1}(0)= \phi (I)=\tau$ is a point. For our $\phi$ derived from Denjoy's example the map $\psi$ is $C^1$-smooth.

Let $ \eta : [-1,1] \rightarrow [0,1]$ be given by the formula $\eta(x) := 1-x^2.$  
Take for $\eta_{\tau}: I \rightarrow [0,\tau]$ the map whose lift is a linear rescaling of $\eta$ to a transformation from $\tilde{I}$ to $ [0,\tau] $.
The following definition gives for each $\sigma \in \{h,v\}$ a continuous $F^{(\sigma)} : X \rightarrow X$. (Here by definition  $ \hat{h}:= v$ and $\hat{v}:=h $.)
\begin{center}
\(
 F^{(\sigma)}(x) := \left\{ \begin{array}{c}
     \psi(\underline{x})^{(\sigma)} \  \mbox{if} \  x \in {\bf S}^{(\sigma)} ;\vspace{0.07in} \\ 
 \eta_{\tau}  (\underline{x})^{(\sigma)} \ \mbox{if} \ x  \in I^{(\hat{\sigma})} ;\vspace{0.07in} \\

p(  \underline{x})^{(\hat{\sigma})} \  \mbox{if} \  x \in S^{(\hat{\sigma})}  \setminus I^{(\hat{\sigma})}.
\end{array} \right.
 \)
\end{center}

\bigskip

\begin{rmk} \hspace{-0.09in}{\bf .} 
The definition is valid for any degree-one mapping $\phi$ which is 1-1 everywhere except a symmetric plateau around zero. Moreover, $F^{(\sigma)} $ depends continuously on $\phi$ in the $C^0$-topology.
\end{rmk}

\bigskip

 The map $F$ is defined as a composition of  $F^{(v)}$ and $ F^{(h)}$,
\[ 
F :=  F^{(v)} \circ F^{(h)}. 
\]

It is easy to homotope  $F^{(h)}$ and $F^{(v)}$ to the identity.
In view of Remark 1 one can do this by coming up with an appropriate homotopy connecting $\phi$ to the identity. Say  $\phi_t :=\psi_t \circ p_t $  where $ \psi_t := id + t \cdot (\psi - id) $ and $p_t$ collapses $t \cdot I$ to zero being affine on the complement of $t \cdot I$.
By lifting homotopies connecting  $F^{(h)}$ and $F^{(v)}$ to the identity to homotopies also terminating at the identity, we get lifts of our mappings that we denote by $\tilde{F}^{(h)}$ and $ \tilde{F}^{(v)}$, respectively. Both $\tilde{F}^{(v)}$ and $\tilde{F}^{(h)}$ belong to $\Hu (X)$. So does their composition  $ \tilde{F}:= \tilde{F}^{(v)} \circ \tilde{F}^{(h)}$.

\begin{prop} \hspace{-0.09in}{\bf .}
 For $\tilde{F}$ defined as above we have 
\[ \rho(\tilde{F}) = \Omega_\rho. \]
\end{prop}

It is convenient to notice that in view of part (iii) of Proposition 1,
Proposition 3 reduces to the following.
\bigskip

\begin{prop} \hspace{-0.09in}{\bf .}  
\begin{itemize}
\item[(i)]  If for $ \tilde{x} \in \tilde{X}$ the limit $\rho(\tilde{F}, \tilde{x} ) = \lim_{n \rightarrow \infty} \frac{1}{n} ( \tilde{F}^n(\tilde{x}) - \tilde{x}) $ exists, then $ \rho (\tilde{F}, \tilde{x}) \in  \Omega_\rho$.
\item[(ii)] For any $ m,n \in {\bf N} \ with \  \alpha_m, \alpha_n < \rho$, there is $\tilde{x} \in \tilde{X}$ with $ \rho (\tilde{F}, \tilde{x}) = {\bf \rho}_{m,n} $.
\end{itemize}
\end{prop}

\bigskip

The remainder of this section is devoted to the proof of Proposition 4. 
\medskip

{\bf Proof of Proposition 4.}
Since $F= F^{(v)} \circ F^{(h)} $ we will find it convenient to regard as the (forward) orbit of $x \in X$ the sequence $ x_0, \ x_{1/2}, \ x_{1}, \ x_{3/2}, \ x_{2},...$, with  
$x_{0}:=x, \ x_{1/2}:=F^{(h)}(x), \ x_{1}:= F(x), \ x_{3/2}:= F^{(h)} \circ F (x)$ ... .

Looking at the definition of the map $F^{(h)}$ we see that it preserves $\Sh$.
The circle  $\Sv$ is not forward invariant, however the only part of $\Sv$ that is folded into $\Sh$ consists of  $I^{(v)}$ which (by definition) is a symmetric neighborhood of the origin. For $F^{(v)}$ we have analogous situation (with $v$ and $h$ switched). 

There are two types of (forward) orbits :

\noindent {\bit eventually free}, i.e. staying in the complement of $I^{(h)} \cup I^{(v)}$ for all sufficiently high indices; 

\noindent  {\bit interacting}, i.e. having returns to $I^{(h)} \cup I^{(v)}$ for arbitrarily large indices. (In fact returns to  $I^{(h)}$ and $I^{(v)}$ must alternate.) 

To understand what happens in each of these two cases and develop useful notation we draw the following diagrams associated with the orbit of $x_0 \in  X$.

\bigskip
{\sl Case of free $x_0 \in \Sh$:} \\

\bigskip
\begin{equation}
 x_{0}^{(h)} \superarrow{ F^{(h)} }{ \psi} x_{1/2}^{(h)} 
\underbrace{\superarrow{F^{(v)}}{p}  x_{1}^{(h)} \superarrow{F^{(h)}}{\psi} }_{\phi} 
x_{3/2}^{(h)} 
\underbrace{ \superarrow{F^{(v)}}{p} x_{2}^{(h)} \superarrow{F^{(h)}}{\psi} }_{\phi} 
x_{5/2}^{(h)}
 \underbrace{ \superarrow{F^{(v)}}{p} x_{3}^{(h)} \superarrow{F^{(h)}}{\psi} }_{\phi}...  \label{diag1}
\end{equation}
\bigskip
\bigskip

{\sl Case of free $x_0 \in \Sv$:} \\

\bigskip

\begin{equation}
 x_{0}^{(v)} 
\underbrace{\superarrow{F^{(h)}}{p}  x_{1/2}^{(v)} \superarrow{F^{(v)}}{\psi} }_{\phi} 
x_{1}^{(v)} 
\underbrace{ \superarrow{F^{(h)}}{p} x_{3/2}^{(v)} \superarrow{F^{(v)}}{\psi} }_{\phi} 
x_{2}^{(v)}
 \underbrace{ \superarrow{F^{(h)}}{p} x_{5/2}^{(v)} \superarrow{F^{(v)}}{\psi} }_{\phi}...
\label{diag2}
\end{equation}
\bigskip
\bigskip

\nopagebreak {\sl Case of interacting  $x_0 \in I^{(v)}$ for which $x_{m+1/2} :=F^{(h)} \circ F^m (x_0)$ sits in $I^{(h)}$ and  $x_{m+n+1} = F^{m+n+1}(x_0)$ is the first return to $I^{(v)}$:} \\
\nopagebreak[3]
\nopagebreak \bigskip
\nopagebreak \bigskip
\nopagebreak 
\begin{eqnarray}
x_{0}^{(v)} \superarrow{ F^{(h)}}{ \eta_\tau} x_{1/2}^{(h)} 
\underbrace{\superarrow{F^{(v)}}{p}  x_{1}^{(h)} \superarrow{F^{(h)}}{\psi} }_{\phi} ...
\underbrace{... \superarrow{F^{(h)}}{\psi} }_{\phi}
x_{m+1/2}^{(h)} \superarrow{F^{(v)}}{\eta_\tau} x_{m+1}^{(v)}
\underbrace{ \superarrow{F^{(h)}}{p} x_{m+3/2}^{(v)} \superarrow{F^{(h)}}{\psi} }_{\phi}   \nonumber \\
x_{m+2}^{(v)}
 \underbrace{ \superarrow{F^{(h)}}{p} ...}_{\phi} ... \underbrace{...  \superarrow{F^{(v)}}{\psi} }_{\phi} x_{m+n+1}^{(v)} \longrightarrow ...
\label{diag3}
\end{eqnarray}

\bigskip
\bigskip
\noindent
The diagrams indicate which of $F^{(h)}$ or $F^{(v)}$ acts and what the corresponding transformation of {\bit circle coordinate} $\underline{x}$  is. For example,
 \[x_{m}^{(h)} \superarrow{F^{(h)}}{\psi} x_{m+1/2}^{(h)} \]
tells us that $x_{m+1/2} = F^{(h)} (x_{m}) $ and $ \underline{x_{m+1/2}} = \psi ( \underline{x_m})$. It also indicates that both $x_{m}$ and $x_{m+1/2}$ belong to $\Sh$.
(Ambiguity arising for $0=0^{(h)}=0^{(v)}$ will cause no confusion in our further discussion so it is ignored.)

We see that if $x_k$ stays on one of $\Ss$, $\sigma \in \{h,v\}$, for $k = 0, \ \frac{1}{2}, \ 1, \  \frac{3}{2}, \ 2, ...$, then the projected coordinate $\underline{x_k} $ on ${\bf S}^1$ evolves under the iterates of $\phi = \psi \circ p$ (i.e. $\underline{x_{k+1}} = \phi( \underline{x_k}$)). When $x_k$ hits $I^{(\sigma)}, \ \sigma \in \{h,v\},$ it makes a transition to the other circle ${\bf S^{(\hat{\sigma})}}$ through the folding action of $\eta_\tau$. 
Given $m,n \in \N$, all points $x_0 \in {\bf S}^{(v)}$ behaving according to diagram (\ref{diag3}) form a set (perhaps empty) that will be denoted $K^{(v)}_{m,n}$. That is, for any $m,n \in \{0,1,2... \}$, we set 
\[
 K_{m,n}^{(v)} := \{ x_0 \in I^{(v)} \ : \ x_{m+n+1} \ \mbox{is the first return to $I^{(v)}$ and} \ x_{m+1/2} \in I^{(h)} \}.
\]

Denote by $h: {\bf S}^1 \rightarrow {\bf S}^1 $ the unique semiconjugacy of $\phi$ to rigid rotation by the angle $\rho$ such that $h(0) = 0$. We define $H : X \rightarrow X $ by $H(x^{(\sigma)}) := h(x)^{(\sigma)}, \ \sigma \in \{h,v\}$. Notice that $H(I^{(h)} \cup I^{(v)}$ ) = (0,0), so in particular, $H$ is a well defined continuous map.
For $\tilde{H} : \tilde{X} \rightarrow \tilde{X} \subset {\bf R}^2 $ we take the lift of $H$ preserving $(0,0)$.

\begin{clm} \hspace{-0.09in}{\bf .} If $\tilde{x} \in \tilde{X}$ is as in (i) of Proposition 4 and $x = \pi(\tilde{x})$, then 
\begin{itemize}
\item[(i)] for $x$ eventually free, $\rho ( \tilde{F}, \tilde{x}) \in \{ (0,\rho),(\rho,0) \}$;
\item[(ii)] for $x$ interacting, $ \rho ( \tilde{F}, \tilde{x}) $ belongs to $ Conv \left\{ \frac{\tilde{H} \circ \tilde{F}^{m+n+1} (\tilde{x}) - \tilde{H} (\tilde{x})}{m+n+1} \ : \ x=\pi (\tilde{x}) \right. $ $\left. \in K_{m,n}^{(v)}, \ m,n \in \N \right\}$.
\end{itemize}
\end{clm}

\bigskip

{\bf Proof of Claim 1.}
We start with a standard observation that for any  $\tilde{H} \in \Hu(X)$ the limits of sequences $\frac{1}{n}(\tilde{F}^n (\tilde{x}) - \tilde{x})$ and $ \frac{1}{n}(\tilde{H} \circ \tilde{F}^n (\tilde{x}) - \tilde{H} (\tilde{x}))$ coincide. Indeed, since $\tilde{H} (\tilde{y}) - \tilde{y} $ is a ${\bf Z}^2$-periodic function of $ \tilde{y}$, we have $\| \tilde{H} \circ \tilde{F}^n (\tilde{x}) - \tilde{F}^n (\tilde{x}) \| \leq \sup_{\tilde{y} \in {\bf R}^2} \| \tilde{H} (\tilde{y}) - \tilde{y} \| < + \infty.$

If $x$ is eventually free, after skipping a finite number of iterates, we are in the situation of one of the two first diagrams.
In the case of diagram (1) we write  $ \tilde{H} \circ \tilde{F}^n (\tilde{x}) - \tilde{H} (\tilde{x}) = \tilde{H} (\tilde{x}_{1/2}) - \tilde{H} (\tilde{x}_0) + \sum_{k=1}^{n-1}(\tilde{H}(\tilde{x}_{k+1/2})-\tilde{H}(\tilde{x}_{k-1/2})) + \tilde{H} (\tilde{x}_{n}) - \tilde{H} (\tilde{x}_{n-1/2}) $.
Since  $\tilde{H} (\tilde{x}_{k+1/2}) - \tilde{H} (\tilde{x}_{k-1/2}) = (\rho,0)$ we get $\rho (\tilde{F},x) =(\rho,0)$.  
Similarly in the case of diagram (\ref{diag2}) we have  $\tilde{H} \circ \tilde{F}^n (\tilde{x}) - \tilde{H} (\tilde{x}) = \sum_{k=1}^{n}(\tilde{H}(\tilde{x}_{k})-\tilde{H}(\tilde{x}_{k-1})) = \sum_{k=1}^{n} (0, \rho)$, so $ \rho (\tilde{F},x) =(0,\rho)$.

For an interacting orbit we may assume that it starts in $I^{(v)}$. Then we split it into segments between consecutive returns to $I^{(v)}$. Each such a segment (after obvious shift of indices) looks like that in diagram (\ref{diag3}). In particular, it has an associated $(m,n) \in {\bf N}^2$ such that $x_{m+1/2} \in I^{(h)}$ and $m+n+1$  is the total length of the segment. Now if we take $N \in {\bf N}$ such that $x_0,..., x_N$ is a certain concatenation of segments as above, then $\frac{1}{N}(\tilde{H} \circ \tilde{F}^N (\tilde{x}) - \tilde{H} (\tilde{x}))$ is a convex combination of average displacements in each segment.
$\Box$
\bigskip

Let us now calculate the average displacement on the sets $K_{m,n}^{(v)}$. 
\medskip
\begin{clm} \hspace{-0.09in}{\bf .} If $\pi(\tilde{x}) \in K_{m,n}^{(v)}$ for some $m,n \in {\bf N}$, then
\[
 \frac{\tilde{H} \circ \tilde{F}^{m+n+1} (\tilde{x}) - \tilde{H} (\tilde{x})}{m+n+1} = {\bf \rho}_{m,n} .
\]
\end{clm}

\bigskip

{\bf  Proof of Claim 2.}
For $\tilde{x}_0$ with $\pi (\tilde{x}_0 ) = x_0 \in K_{m,n}^{(v)} $, we have  \\
\(
 \tilde{H}( \tilde{x}_{m+n+1}) - \tilde{H}(\tilde{x}_{0}) =  (\tilde{H}( \tilde{x}_{1/2}) - \tilde{H}(\tilde{x}_{0}) ) + \sum_{k=1}^{m} (\tilde{H}( \tilde{x}_{k+1/2}) - \tilde{H}(\tilde{x}_{k-1/2}) )+ (\tilde{H}( \tilde{x}_{m+1}) - \tilde{H}(\tilde{x}_{m+1/2}) +
\sum_{l=1}^{n} (\tilde{H}( \tilde{x}_{m+l+1}) - \tilde{H}(\tilde{x}_{m+l}) ) .
\)

Observe that  $H(x_0) = H(x_{m+1/2}) = H(x_{m+n+1})=(0,0)$. Also $H( \phi^{-r} (I) ^{(\sigma)} ) = (-r \rho \ (mod \ {\bf Z}))^{(\sigma)}$ for  $r \in {\bf N}$, while $x_{k+1/2} \in (\phi^{-m+k}(I))^{(h)}$ for $k=0,...,m$ and  $ x_{m+l+1} \in (\phi^{-n+l}(I))^{(v)}$ for $l=0,...,n$. Thus $ \tilde{H}( \tilde{x}_{m+n+1}) - \tilde{H}(\tilde{x}_{0}) = (-m \rho - \lfloor -m\rho \rfloor) \cdot
(1,0) + m \rho \cdot (1,0) +  (-n\rho - \lfloor -n\rho \rfloor) \cdot (0,1) + n \rho \cdot (0,1)=(\lceil m\rho \rceil,\lceil n \rho \rceil) = (m+n+1) \cdot \rho_{m,n}$. 
$ \Box$

\bigskip

\begin{lem} \hspace{-0.09in}{\bf .}
The set $K_{m,n}^{(v)} $ is nonempty if and only if $\alpha_m, \alpha_n < \rho$. Moreover, if $K_{m,n}^{(v)} $ is nonempty, then it is a disjoint union of four subarcs of $I^{(v)}$, each mapped by $ F^{m+n+1}$ homeomorphically onto $I^{(v)}$.
\end{lem}

\bigskip

{\bf Proof of Lemma 1.} The arcs $\{ \phi ^{-j}(I) ^{(\sigma )} \}_{j \in \N} $ are pairwise disjoint and ordered on $S^{(\sigma )}$ in the same manner as the corresponding  backward orbit of rigid rotation by $\rho $.
 Also $\eta _\tau (I) = [0, \phi (I)]$. 
In this way the set of all points in $I$ which get folded by $\eta _\tau$ into $\phi ^{-m}(I)$ is either empty or consists of a two arcs, each mapped by $\eta _\tau$ homeomorphically onto $\phi^{-m}(I)$.
The set is nonempty exactly when  $-m\rho \in (0,\rho ) \ (mod \ \Z)$, i.e. $-m \rho - \lfloor -m \rho  \rfloor < \rho $ or equivalently $\alpha _m = \ceil { m \rho } - m \rho < \rho $.

Now look at diagram (3). 
For a point $x_0 \in I^{(v)}$ to have $x_{m+1/2}$ as the first iterate
hitting $I^{(h)}$  means that $\eta _\tau (\underline{x_0}) \in \phi ^{-m}(I)$. In view of the preceding remarks such points $x_0$ exist only if $\alpha_m < \rho$, and then they form two subarcs of $I^{(v)}$, call them $K_1$ and $K_2$.  
Note that under  $F^{m+1/2} = F^{m} \circ F^{(h)}$ (which is essentially $\phi^m \circ \eta_\tau $ ), each of them maps homeomorphically onto $I^{(h)}$. 

Analogously, points $x_{m+1/2}$ in $I^{(h)}$ that have $x_{m+n+1}$ as the first iterate hitting $I^{(v)}$ exist if and only if $\alpha_n < \rho$. As before, if $\alpha_n < \rho$, then those points form two subarcs of $I^{(h)}$, call them $L_1$ and $L_2$. Under  $F^{n} \circ F^{(v)}$ (which is essentially $\phi^n \circ \eta_\tau $ ), each of them maps homeomorphically onto $I^{(v)}$ (see diagram (3)).

Now $K_{m,n}^{(v)}$ is clearly equal to $F^{-(m+1/2)}(L_1 \cup L_2)\cap (K_1\cup K_2)$. This set, if it is nonempty, consists of four subarcs of $I^{(v)}$ each mapped by $F^{m+n+1}=F^{m+1/2} \circ F^{(v)} \circ F^{n}$  homeomorphically onto $I^{(v)}$. The lemma is proved.
$\Box$

\bigskip

We are ready to conclude the proof of Proposition 4.
To prove (ii), note that from Lemma 1 we see that for $m$, $n$ with $ \alpha_m, \ \alpha_n < \rho$ the sets $K_{m,n}^{(v)}$, $F^{-(m+n+1)}(K_{m,n}^{(v)})$, $F^{-2(m+n+1)}(K_{m,n}^{(v)}), ...$ are nested, so the intersection  $ \bigcap_{k=0}^{\infty} F^{-k(m+n+1)}(K_{m,n}^{(v)}) $ is nonempty. Take any  $x$ in this set.

 To prove (i) we have to consider two possibilities. For eventually free $x$ we are done by part (i) of Claim 1 together with the fact that both $(0, \rho)$ and $(\rho,0)$ belong to $\Omega_{\rho}$.
For interacting $x$ we have to combine part (ii) of Claim 1 with Claim 2.
$\Box$

\bigskip
\bigskip
\begin{center}
{\large\bf Section 3: A diffeomorphism of the torus from the map on the one-skeleton.}
\end{center}
\bigskip

Recall that  $ \frac{18}{17} \cdot I$ is a wandering arc under $\phi$ and set  \[ 
W := \phi \left( \frac{18}{17} \cdot I \right).
\]
Shrinking $\frac{18}{17} \cdot I $ we can  get a symmetric arc $ I''$ containing $I$ such that  $ J_{0} := p(I'')$  is a nondegenerated arc and $\eta_\tau ( J_0)$ is contained in W (see Figure 2). Note that $ \psi ( J_0) \subset W$. Put 
\[
I' := \eta_{\tau}^{-1} ( [0, \tau] \setminus J_0)
\]
 and 
\[
 \Delta I := cl(I'' \setminus I').
\]
(See Figure 2.)

\begin{figure}[htp]
\centerline{\psfig{figure=folding1.fig,width=\hsize}}
\centerline{Figure 2.}
\end{figure}

For $ \sigma \in \{h,v\} $ put  
\[
U^{(\sigma)} := \{ x \in {\bf T}^2 \ : \ dist(x,\Ss) < 1/4 \cdot  length ( J_0)\}; 
\]
\[ 
 U:= U^{(h)} \cup U^{(v)}.
\] 
Also set 
\[
 C :=  \{  (x,y) \in U \ : \ |x|, |y| < 1/2 \cdot length (J_0) \}.
\]

\noindent
There is a retraction  $ r : U \setminus C \rightarrow X \setminus C $ given by 
\[
 r(x,y):= \left\{ \begin{array}{c} (x,0) \mbox{ for } (x,y) \in U^{(h)}; \\
(0,y) \mbox{ for } (x,y) \in  U^{(v)}. \end {array} \right. 
\]

Theorem 1 is a consequence of the following proposition which is the main result of this section.

\begin{prop} \hspace{-0.09in}{\bf .}
There exists a $C^1$-embedding $G : cl(U) \rightarrow U$ such that $\rho(G) = \Omega_{\rho}$.
\end{prop}

\bigskip

Let us first see how Theorem 1 follows from this proposition.

{\bf Proof of Theorem 1.}
It is a routine to  extend the map $G$ acting on $U \subset {\bf T}^2$ to a diffeomorphism on the whole ${\bf T}^2$ by putting a single source outside $U$ repelling all other points towards $ \bigcap_{n \geq 0} G^{n} (U) $ (i.e. we require that  $ {\bf T}^2 \setminus \bigcap_{n \geq 0} G^{n}(U) $ is a basin of attraction under $G^{-1}$).
For this extended $G$ we have $\rho(G) = Conv \{(0,0), \Omega_\rho \} = \Lambda_\rho$.
$\Box$

\bigskip
  
The rest of this section is devoted to the proof of Proposition 5. 

\medskip

{\bf Proof of Proposition 5.}
The construction of the $C^1$-embedding $G : cl(U) \rightarrow U$ 
is done in steps.

\bigskip

{\bf Step 1.} 
As suggested by Figure 3, we  modify  $F^{(h)} $ on $ (I'')^{(v)} $ to obtain $ F^{(h)}_1 : X \rightarrow  U$ a $C^1$-smooth embedding satisfying the conditions listed below.

\begin{figure}[htp]
\centerline{\psfig{figure=folding2.fig,width=\hsize}}
\centerline{Figure 3.}
\end{figure}

\begin{itemize}
\item[(A)] $ (F^{(h)}_1)^{-1}(cl(C)) = \Delta I^{(v)} \cup \psi^{-1}(J_0)^{(h)}$;
\item[(B)] $F^{(h)}_1 ( cl(C)) \subset r^{-1}( W ^{(h)})$;
\item[(C)] if $x \nin C \cup (F^{(h)}_1)^{-1}(C)$, then $ r \circ F^{(h)}_1 (x) = F^{(h)}(x). $ 
\end{itemize}

\bigskip

{\bf Step 2.}
We extend $F^{(h)}_1$ to a $C^1$-smooth embedding  $G^{(h)} :
cl(U) \rightarrow U $ with the following analogues of conditions (A), (B), (C) satisfied :
\begin{itemize}
\item[(A$^\prime$)] $ (G^{(h)})^{-1}(cl(C)) = r^{-1}( \Delta I^{(v)} \cup \psi^{-1}(J_0)^{(h)})$;
\item[(B$^\prime$)] $G^{(h)} ( cl(C)) \subset r^{-1}( W ^{(h)})$;
\item[(C$^\prime$)] if $x \nin C \cup (G^{(h)})^{-1}(C)$, then $ r \circ G^{(h)} (x) = F^{(h)} \circ r(x).$
\end{itemize}

\bigskip
See Figure 4  for the construction of $G^{(h)}$. In the analogous way we obtain $G^{(v)}$ from $F^{(v)}$.

\begin{figure}[htp]
\centerline{\psfig{figure=folding3.fig,width=\hsize}}
\centerline{Figure 4.}
\end{figure}

\bigskip

{\bf Step 3.} Finally we define $G : cl(U) \rightarrow U$ by 
\[
G := G^{(v)} \circ G^{(h)} .
\]

\bigskip

Let us now concentrate on these aspects of the dynamics of 
$G$ that are critical for the calculation of the rotation set $\rho(G)$.
First let us  see how alternating applications of $G^{(v)}$ and $ G^{(h)} $ act on $r^{-1}(W^{(h)})$ and $r^{-1}(W^{(v)})$.
We expect that these domains are wandering under G and move ``freely'' governed by the action of $\phi$ on the appropriate circle. 
We formalize this in the following lemma.
(We deal with the case of $r^{-1}(W^{(h)})$. For $r^{-1}(W^{(v)})$ there is an analogous statement.)

\bigskip
\begin{lem} \hspace{-0.09in}{\bf .} 
For any $M = 0,1,2...$  
\begin{itemize}
\item[(a)]  if M is even and $M=2 \cdot N$, then the set $\overbrace{ G^{(h)} \circ G^{(v)} \circ ... \circ G^{(h)} \circ G^{(v)} }^{2N \ compositions} ( r^{-1} (W^{(h)})) $  is contained in $r^{-1}( \phi^N(W)^{(h)})$ which is disjoint from $C \cup (G^{(v)})^{-1}(C) $; 
\vspace{0.06in}
\item[(b)] if M is odd and  $M=2 \cdot N +1$, then   $\overbrace{ G^{(v)} \circ G^{(h)} \circ ... \circ G^{(h)} \circ G^{(v)} }^{2N+1 \ compositions } ( r^{-1} (W^{(h)})) $  is contained in $r^{-1}(p \circ \phi^N(W)^{(h)})$ which is disjoint from $C \cup (G^{(h)})^{-1}(C) $.
\end{itemize}
\end{lem}

\bigskip

{\bf Proof of Lemma 2.}
We proceed by induction on M. 

For $M=0$ our claim follows from (A$^\prime$) for $G^{(v)} $ and the fact that $W \cap (J_0 \cup \Delta I) = \emptyset$. 

Now we will describe the induction step.  

First we deal with part (b). 
By the inductive hypothesis and  (C$^\prime$) for  $ G^{(v)}$  we have  $ r \circ \overbrace{ G^{(v)} \circ ... \circ G^{(v)} }^{2N+1} ( r^{-1}(W^{(h)})) \subset 
r \circ G^{(v)} \circ r^{-1}( \phi^N (W) ^{(h)}) = r \circ r^{-1} \circ F^{(v)}( \phi^N (W) ^{(h)}) =  p \circ \phi^N(W)^{(h)} $. 
Thus  $ \overbrace{ G^{(v)} \circ ... \circ G^{(v)} }^{2N+1} ( r^{-1}(W^{(h)})) \subset
 r^{-1}(p \circ \phi^ N ( W )^ {(h)}) $.

Now  $r^{-1} (p \circ \phi^N(W)^{(h)}) \cap C = \emptyset $ follows from $p \circ \phi ^N(W) \cap J_0 = \emptyset$ which is equivalent to $ \psi \circ p \circ \phi^N (W) \cap \psi(J_0) = \phi^{N+1}(W) \cap \psi(J_0) = \emptyset $. Since  $ \psi( J_0) \subset W $ and W is wandering this is true. 

On the other hand, $ r^{-1 } ( p \circ  \phi^{N} ( W ^ {(h)}))  \cap ( G^{(h)})^{-1}(C) = \emptyset$ in view of (A$^\prime$) follows from  $ p \circ \phi^N(W) \cap \psi^{-1}( J_0) = \emptyset $ or equivalently $\phi^{N+1}(W) \cap J_0 = \emptyset $, once again by the wandering property of $W$. 

Let us now consider part (a).
By the inductive hypothesis  and  (C$^\prime$) for  $ G^{(h)}$ we have $   r \circ \overbrace{ G^{(h)} \circ ... \circ G^{(v)} }^{2N} ( r^{-1}(W^{(h)}) = r \circ G^{(h)} \circ r^{-1} ( p \circ \phi^{N-1}(W)^{(h)}) = r \circ r^{-1} \circ F^{(h)} (p \circ  \phi ^{N-1}(W)^{(h)}) =  ( \psi( p \circ \phi ^{N-1}(W))^{(h)} = \phi^N (W) ^{(h)}$.\\
Now  $ r^{-1}( \phi^N (W) ^{(h)}) \cap C = \emptyset $ follows from $\phi^N(W) \cap J_0 = \emptyset.$ 
On the other hand, using the version of (A$^\prime$) for $G^{(v)}$ we see that $r^{-1}( \phi^N (W) ^{(h)}) \cap (G^{(v)})^{-1}(C) = \emptyset$ follows from $\phi ^N (W) \cap \Delta I = \emptyset$. $\Box$

\bigskip

Now we are going to verify that the region $C$ under iterations of $G$ stays in the wandering domains of Lemma 2.
This assures that all the rotationally nontrivial dynamics is carried on a certain forward invariant set $\Sigma$ on which $r$ is a well defined semiconjugacy between $G$ and $F$. The lemma below gives precise formulations.

\begin{lem}  \hspace{-0.09in}{\bf .}
If  $Crit := \bigcup_{n=0}^{\infty} C_{-n/2} = \bigcup_{n=0}^{\infty} G^{-n}( C \cup (G^{(h)}) ^{-1}(C)) $ and $ \Sigma :=U \setminus Crit$, then
\begin{itemize}
\item[(i)]  $G(\Sigma) \subset \Sigma$ ;
\item[(ii)]  For any $x \in U$ there is  $N \in \{0,1,2...\}$ such that $ x_N := G^N(x)  \in \Sigma$. Moreover, if $x \nin \Sigma$ then there exists $N \in {\bf N}$ such that either $x_{N+k} \in  r^{-1}(p \circ \phi^k (W)^{(h)})$ for all $k \in {\bf N}$ or  $x_{N+k} \in  r^{-1}( \phi^k (W)^{(v)})$ for all $k \in {\bf N}$; 
\item[(iii)] The following diagram commutes : 
\[
\begin{array}{ccc}
\Sigma & \stackrel{G}{\longrightarrow} & \Sigma \\
r \downarrow &   & \downarrow r \\
r( \Sigma ) & \stackrel{F}{\longrightarrow} & r( \Sigma) \\
\end{array}
\]
\end{itemize}
\end{lem}

Note that the nonwandering set of $G$ is contained in $\Sigma$. Also the nonwandering set of $F$ sits in $r(\Sigma)$.

\bigskip

{\bf Proof of Lemma 3}. 

\noindent
(ii)   Suppose that $ x_0 \nin \Sigma$ . Then for some $N \in \{1,2, ... \}$  we have $ x_{N-1} \in C \cup (G^{(h)})^{-1}(C)$.
 
If $x_{N-1} \in C$, then by (B$^\prime$) we have $x_{N-1/2} = G^{(h)} (x_{N-1}) \in r^{-1}(W^{(h)}) $ and from  Lemma 2 we conclude that 
for every $k \geq 0 $  \[  x_{N+k} = \overbrace{ G^{(v)} \circ G^{(h)} \circ ... \circ G^{(h)} \circ G^{(v)} }^{2k+1 \ compositions} (x_{N-1/2}) \] is contained in $r^{-1}( p \circ \phi^k(W) ^{(h)}) $. This set is disjoint from $C \cup (G^{(h)})^{-1}(C)$ so $x_N \nin Crit$. 

If $x_{N-1} \in (G^{(h)})^{-1}(C)$, then $x_{N-1/2} = (G^{(h)})(x_{N-1}) \in C $ and, by (B$^\prime$) for $G^{(v)}$, we have $x_N = G^{(v)} (x_{N-1/2}) \in r^{-1}(W^{(v)})$. Applying Lemma 2 with the role of $v$ and $h$ interchanged we see that 
for every $ k \geq 0$  \[ x_{N+k} = \overbrace{ G^{(v)} \circ G^{(h)} \circ ... \circ G^{(v)} \circ G^{(h)}}^{2k \ compositions}  (x_N) \]  is contained in $r^{-1}( \phi^k (W)^{(v)}) $.  This set is disjoint from $C \cup(G^{(h)})^{-1}(C)$ and so $x_N \nin Crit$.

\bigskip
(iii)     Suppose that $x_0 \nin Crit$, then $x_0 \nin C \cup (G^{(h)})^{-1}(C)$ and $x_{1/2} = G^{(h)} (x_0)  \nin C \cup (G^{(v)})^{-1}(C)$. Indeed, otherwise $x_0 \in C \cup (G^{(v)})^{-1}(C) \cup (G^{(v)} \circ G^{(h)})^{-1}(C) \subset Crit$. Thus we may apply (C$^\prime$) with $x=x_0$ and (C$^\prime$) for $G^{(v)}$ with $x=x_{1/2}$ to get  $r \circ G^{(v)} \circ G^{(h)} (x) = F^{(v)} \circ F^{(h)} (r (x))$. $\Box$

\bigskip

To conclude the proof of Proposition 5 we need to show that $ \rho(F) = \rho(G) = \Omega_\rho$.
Observe that Lemma 3 part (ii) tells us that points $x \nin \Sigma$ contribute only $(0,\rho)$ or $(\rho,0)$ to the rotation set of $G$.
The analogous statement is also true for $F$, that is, the contribution of points which are not in $r(\Sigma)$ to the rotation set of $F$ is either $(\rho,0)$ or $(0, \rho)$. This contribution is due to free orbits on $\Sh$ or $\Sv$ correspondingly. Moreover, using  the diagram in (iii) of Lemma 3 lifted to the universal cover we see that for any $x \in \Sigma$,  $\rho( \tilde{G},x) = \rho(\tilde{F}, r(x))$.  Indeed, we have $ \|\tilde{G}^{n} (\tilde{x}) -
 \tilde{F}^{n} \circ \tilde{r} (\tilde{x}) \| = \| \tilde{G}^{n} (\tilde{x}) - \tilde{r} \circ \tilde{G}^{n} (\tilde{x}) \| \leq \sup_{ \tilde{y} \in \tilde{\Sigma}} \| \tilde{r}( \tilde{y})- \tilde{y} \| < \infty$.
$\Box$

\bigskip
\bigskip
\begin{center}
{\large\bf Section 4: Deriving invertible dynamics on the torus from noninvertible dynamics on a skeleton.}
\end{center}
\bigskip

In this section we want to indicate that, as far as one is ready to give up on  smoothness requirements, there is a general method essentially replacing considerations of Section 3. The construction is fairly robust and intuitive, so we assume a very informal style to avoid blurring the essence with a cloud of details.

Our approach is a variation of the method used by Barge and Martin (\cite{BM}) to prove that inverse limits of interval transformations can be realized as attractors for homeomorphisms of the plane. 
The key fact is the following theorem which is an easy corollary of Morton Brown's results in  \cite{Brown}. Let us recall that a continuous map $f$ of a compact metric space $Y$ is a  {\bit nearhomeomorphism} if and only if there exists a sequence of homeomorphisms $f_n : Y \rightarrow Y$ converging uniformly to $f$. 

\begin{thm} \hspace{-0.09in}{\bf .}
If $\, Y$ is a compact metric space, then any $nearhomeomorphism$ $f : Y \rightarrow Y$ is a factor of a homeomorphism $\hat{f}$ of $\, Y$, i.e. there is a continuous onto map $h :Y \rightarrow Y$ such that $ h \circ \hat{f} = f \circ h$. 
\end{thm}

The paper of Brown is clear and puts the above theorem in a perspective of general facts about inverse limits. However, we feel it will be  beneficial to the reader if we present here a self-contained proof of the result.

\bigskip

{\bf Proof of Theorem 2.}  Let $d$ be the metric on $Y$. Given any two maps of Y, say  $f$ and $g$, use $d(f,g)$ to denote their uniform distance equal by definition to $\sup\{d(f(x),g(x)) : x \in Y\}$. Suppose that a sequence of homeomorphisms $(f_n)_{n \in \N}$ converges uniformly to $f$. Compositions $f_1 \circ ...\circ f_n$ and $f_n^{-1} \circ ... \circ f_1^{-1}$ are denoted by $f_{1,n}$ and $f_{n,1}^{-1}$ respectively. 
For any $n \in {\bf N}$  set
\[ \alpha_n := \sup\{ \max \{d(f_{1,n-1}(x),f_{1,n-1}(y)), \ d(f_{1,n-1} \circ f(x) ,f_{1,n-1} \circ f (y)) \} \] \[ : \ d(x,y) < d(f_n,f)\}, \]
\[ \beta_n:= \sup\{ \max \{d(f_{1,n-1}(x),f_{1,n-1}(y)), \ d(f_{1,n-1} \circ f (x), f_{1,n-1}\circ f (y)) \} \] \[  : \ d(x,y) < d(f_{n+1},f)\},
\]
\[ \gamma_n := sup \{d(f^{n-1}(x),f^{n-1}(y)) : d(x,y)<d(f,f_n) \}, \]
\[ \epsilon_n := 3 \cdot ( \alpha_n+ \beta_n+\gamma_n). \]
Notice that skipping certain elements of the sequence $(f_n)_{n=1}^{\infty}$ one can make the sequence $(\epsilon_n)_{n=1}^{\infty}$ converging to zero as fast as we wish. Indeed, one can proceed inductively as follows. Suppose that we have already chosen $f_1,...,f_{n}$ so that $\alpha_1,...,\alpha_n$, $\beta_1,...,\beta_{n-1}$ and $\gamma_1,...,\gamma_n$ are as small as we wished. By the uniform continuity of the functions $f_{1,n}$, $f_{1,n-1}$, $f_{1,n} \circ f$, $f_{1,n-1}  \circ f$ and $f^n$, by picking $f_{n+1}$ sufficiently close to $f$ we get $\alpha_{n+1}$, $\beta_n$ and $\gamma_{n+1}$ as small as we wish. This ends the induction step. In this way we may  assume that $\sum \epsilon_n < \infty$.

We claim that the sequences consisting of the following maps of X are  uniformly convergent :
\[ \hat{f}_n := f_{1,n+1} \circ f_{n,1}^{-1}, \]
\[ \hat{g}_n := f_{1,n} \circ f_{n+1,1}^{-1}, \]
\[ h_n := f^{n} \circ f_{n,1}^{-1}. \]

This claim already implies the theorem. Indeed, we have $\hat{g}_n \circ \hat{f}_n = \hat{f}_n \circ \hat{g}_n = id_X$ and $h_{n+1} \circ \hat{f}_{n} = f \circ h_{n}$. Consequently, the limits, denoted by $\hat{f}$, $\hat{g}$, $h$ correspondingly, satisfy $\hat{f} \circ \hat{g} = \hat{g} \circ \hat{f} = id_{X}$ and $h \circ \hat{f} = f \circ h$.

Convergence of the sequences is an immediate consequence of the assumption about the convergence of the series of $\epsilon_{n}$'s and  the following estimates valid for any $x \in X$.
Setting $y:= f_{n,1}^{-1}(x)$ we see that
\[ d(\hat{f}_n(x), \hat{f}_{n-1}(x)) = \]
 \[
 d( f_{1,n+1} \circ f_{n,1}^{-1}(x) , f_{1,n} \circ f_{n-1,1}^{-1}(x)) =  d( f_{1,n-1} \circ f_n \circ f_{n+1}(y), f_{1,n-1} \circ f_n \circ f_n (y)) \leq \] \[
d( f_{1,n-1} \circ f \circ f_{n+1}(y), f_{1,n-1} \circ f \circ f (y)) +
d( f_{1,n-1} \circ f \circ f_{n}(y) , f_{1,n-1} \circ f \circ f (y))+ \] \[
d( f_{1,n-1} \circ f_n \circ f_{n+1}(y) , f_{1,n-1} \circ f \circ f_{n+1} (y)) +
d( f_{1,n-1} \circ f_n \circ f_{n}(y), f_{1,n-1} \circ f \circ f_n (y)) \leq \] \[
 \beta_n + \alpha_n + \alpha_n + \alpha_n \leq \epsilon_n. \]
Similarly taking $f_{n+1,1}^{-1}(x)$ for $y$ we get
\[ d(\hat{g}_n(x),\hat{g}_{n-1}(x)) =
 \]
\[ d( f_{1,n} \circ f_{n+1,1}^{-1}, f_{1,n-1} \circ f_{n,1}^{-1}) =
   d( f_{1,n-1} \circ f_{n} (y) , f_{1,n-1} \circ f_{n+1} (y)) \leq
 \]
\[ d( f_{1,n-1} \circ f_n (y), f_{1,n-1} \circ f (y)) +
   d( f_{1,n-1} \circ f_{n+1} (y), f_{1,n-1} \circ f (y)) \leq 
\]
\[ \beta_n+ \alpha_n \leq \epsilon_n.
\]
Finally  using $y:= f_{n,1}^{-1}(x)$ we obtain 
\[
 d( h_{n}(x), h_{n-1}(x)) =
\]
\[
 d( f^{n} \circ f_{n,1}^{-1}(x),  f^{n-1} \circ f_{n-1,1}^{-1}(x)) =
  d( f^{n-1} \circ f (y), f^{n-1} \circ f_n (y)) \leq 
\]
\[ 
\gamma_n \leq \epsilon_n .
\]
 $\Box$

\bigskip

\begin{figure}[htp]
\centerline{\psfig{figure=foldingsup.fig,width=\hsize}}
\centerline{Figure 5.}
\end{figure}

We will now explain how Theorem 2 enables us to find a torus homeomorphism 
with rotation set equal to $\Lambda_{\rho}$. Section 1 and Section 2 are prerequisites for our considerations. In particular, one should look there for definitions. 
For the space $Y$ in Theorem 2 we take the torus ${\bf T}^2$. For the map $f$ we need a homotopic to the identity nearhomeomorphism which is an extension of the map $F: X \rightarrow X$ to the whole torus and has $X=\Sh \cup \Sv$ as a global attractor. Such a map $f$ would have the desired rotation set. Construction of the map $f$ is not difficult. Once again it is convenient to obtain it as a composition of two maps $f^{(h)}$ and $f^{(v)}$. 
We will roughly sketch the construction of $f^{(h)}$ now. To get $f^{(v)}$ interchange the role of $horizontal$ and $vertical$ directions.

First we take an  embedding $g^{(h)}$ of $X$ into  the torus such that its postcomposition with the map collapsing the marked annulus (see Figure 5) onto $\Sh$ along vertical fibers is equal to $F^{(h)}$. Next we extend $g^{(h)}$ to a homeomorphism of $\T$ which we are going to denote with the same letter. We also extend the collapsing map to a mapping $p^{(h)}$ of ${\bf T}^2$. These extensions need to be reasonably chosen since we want $X$ to be a global attractor. For $f^{(h)}$ we take the composition $p^{(h)} \circ g^{(h)}$. This is a nearhomeomorphism because $p^{(h)}$ is one: the maps $p_n^{(h)}$ contracting the annulus along the vertical direction by a factor say $1/n$ are homeomorphisms that approximate it.

From Theorem 2 we get a homeomorphism $\hat{f}$ of which $f$ is a factor. This map has the same rotation set as $f$ - it is  equal to $\Lambda_{\rho}$ as required. 

The last implication hinges on the fact that the factor map is homotopic to the identity and the following simple observation.

\begin{fct} \hspace{-0.09in}{\bf .}
If $f, \hat{f}, r: \T \rightarrow \T$ are continuous, homotopic to the identity and satisfy $f \circ r = r \circ \hat{f}$, and furthermore $F, \hat{F}, R: {\bf R}^2 \rightarrow {\bf R}^2$ are their lifts such that $F \circ R = R \circ \hat{F}$, then $\rho(\hat{F})=\rho(F)$.
\end{fct}

{\bf Proof.} The map $R$ is surjective and, for any $x \in {\bf R}^2$, we have $ \| \hat{F}^{n} (x) - F^{n} \circ R (x) \| = \|\hat{F}^{n}(x) - R \circ \hat{F}^{n} (x) \| \leq \sup_{ y \in \R } \| R(y) - y \| < \infty$. The fact follows now trivially from the definition of the rotation set. 
$\Box$

If $r$ is not homotopic to the identity, it induces a linear map $r_*$ on the first real homology which we identify with the universal cover ${\bf R}^2$.  Then $\rho(F)=r_*(\rho(\hat{F}))$, provided $r_*$ is nonsingular (or at least $r$ is surjective).
To see this one can modify slightly the above proof, or just compose $r$ with the linear torus automorphism covered by $r_*^{-1}$ and use Fact 1.

\medskip

Let us end with a remark on the continuity of the dependence of 
$\hat{f}$ on $f$.  

\begin{rmk} \hspace{-0.09in}{\bf .}
If $f^{t}:Y \rightarrow Y,\ t \in [0,1]$ is a homotopy of nearhomeomorphisms and there exists a sequence of isotopies  $f_n^t : Y \rightarrow Y, \ t \in [0,1]$ converging uniformly in t to $f^t$, then one can assure that the family of maps $\hat{f}^t$ from Theorem 2 is also a homotopy.
\end{rmk}
 
To see why the remark is true, observe that the dependence of maps  $\hat{f}^t$ on $t$ is continuous if only all estimates in the proof of the theorem are uniform in $t$.
This uniformity would follow if we choose the approximating family $f^t_n$ in the proof  so that $\epsilon_n := \sup\{\epsilon^t_n: \ t \in [0,1] \}$ has as before a finite sum  $\sum \epsilon_n < \infty$.  
One achieves this by analogous inductive procedure as that from the beginning of the argument.

\bigskip
\bigskip
{\bf Acknowledgments.} I would like to thank John Milnor for 
asking the question that prompted this work. I am also very much indebted to Phil Boyland whose support and countless suggestions made this paper possible.

\bibliography{references}
\bibliographystyle{amsplain}
\end{document}